\documentclass[reqno]{amsart}
\usepackage{amssymb}

\usepackage{amsthm}
\usepackage{amsmath}
\usepackage[all]{xy}
\usepackage{color}
\usepackage{enumerate}
\usepackage{fancyhdr}

\newtheorem{theorem}{Theorem}[section]
\newtheorem{proposition}[theorem]{Proposition}
\newtheorem{corollary}[theorem]{Corollary}
\newtheorem{lemma}[theorem]{Lemma}
\theoremstyle{definition}

\newtheorem*{Beweis}{Proof}
\newtheorem{definition}[theorem]{Definition}
\newtheorem{example}[theorem]{Example}
\newtheorem{examples}[theorem]{Examples}
\newtheorem{punto}[theorem]{}
\theoremstyle{remark}
\newtheorem{remark}[theorem]{Remark}

\newtheorem{remarks}[theorem]{Remarks}
\swapnumbers
\input{tcilatex}

\begin{document}
\title[Tilting Modules over APD's]{Tilting Modules over Almost Perfect Domains}
\author{Jawad Abuhlail}
\address{Department of Mathematics and Statistics\\
Box 5046, KFUPM, 31261 Dhahran, KSA}
\email{abuhlail@kfupm.edu.sa}
\urladdr{http://faculty.kfupm.edu.sa/math/abuhlail}
\thanks{The authors were supported by KFUPM under Research Grant \# MS/Rings/351.}
\author{Mohammad Jarrar}
\address{Department of Mathematics and Statistics\\
Box 5046, KFUPM, 31261 Dhahran, KSA}
\email{mojarrar@kfupm.edu.sa}
\urladdr{}
\date{}
\maketitle

\begin{abstract}
We provide a complete classification of all \emph{tilting modules} and \emph{%
tilting classes} over almost perfect domains, which generalizes the
classifications of tilting modules and tilting classes over Dedekind and $1$%
-Gorenstein domains. Assuming the APD is Noetherian, a complete
classification of all \emph{cotilting modules} is obtained (as duals of the
tilting ones).
\end{abstract}

\section{Introduction}

Throughout, $R$ is a commutative ring with $1_{R}\neq 0_{R}$ and all $R$%
-modules are unital. With $Z(R)$ we denote the set of zero-divisors of $R$
and set $R^{\times }:=R\backslash Z(R).$ With $Q=(R^{\times })^{-1}R$ we
denote the total ring of quotients of $R$ (the field of quotients, if $R$ is
an integral domain). With $R$\textrm{-}$\mathrm{Mod}$ we denoted the
category of $R$-modules.

Let $M$ be an $R$-module. The \emph{character module} of $M$ is $M^{c}:=%
\mathrm{Hom}_{\mathbb{Z}}(M,\mathbb{Q}/\mathbb{Z}).$ With $\mathrm{Max}(M)$
we denote the (possibly empty) spectrum of maximal $R$-submodules and define
\begin{equation*}
\mathrm{rad}(_R M):=\bigcap\limits_{L\in \mathrm{Max}(M)}L\text{ }\;\text{(}%
=M,\text{ if }\mathrm{Max}(M)=\varnothing \text{).}
\end{equation*}
In particular, $\mathrm{Max}(R)$ is the spectrum of maximal $R$-ideals and $%
J(R):=\mathrm{rad}(_R R)\;$is the Jacobson radical of $R.$ We denote with $%
\mathrm{p.d.}_{R}(M)$ (resp. $\mathrm{i.d.}_{R}(M),$ $\mathrm{w.d.}_{R}(M)$)
the projective (resp. injective, weak or flat) dimension of $_{R}M.$
Moreover, we set
\begin{equation*}
\begin{tabular}{lllllll}
$\mathcal{P}_{n}$ & $:=$ & $\{_{R}M\mid \mathrm{p.d.}_{R}(M)\leq n\};$ &  & $%
\mathcal{P}$ & $:=$ & $\dbigcup\limits_{n=0}^{\infty }\mathcal{P}_{n};$ \\
$\mathcal{I}_{n}$ & $:=$ & $\{_{R}M\mid \mathrm{i.d.}_{R}(M)\leq n\};$ &  & $%
\mathcal{I}$ & $:=$ & $\dbigcup\limits_{n=0}^{\infty }\mathcal{I}_{n};$ \\
$\mathcal{F}_{n}$ & $:=$ & $\{_{R}M\mid \mathrm{w.d.}_{R}(M)\leq n\};$ &  & $%
\mathcal{F}$ & $:=$ & $\dbigcup\limits_{n=0}^{\infty }\mathcal{F}_{n}.$%
\end{tabular}
\end{equation*}
In particular, $\mathcal{PR}:=\mathcal{P}_{0}$ is the class of projective $R$%
-modules, $\mathcal{IN}:=\mathcal{I}_{0}$ is the class of injective $R$%
-modules, and $\mathcal{FL}:=\mathcal{F}_{0}$ is the class of flat\emph{\ }$%
R $-modules. The class of torsion-free $R$-modules will be denoted with $%
\mathcal{TF}.$ For a multiplicative subset $S\subseteq R^{\times },$ the
class of $S$\emph{-divisible} $R$-modules is
\begin{equation*}
\mathcal{D}_{S}:=\{_{R}M\mid sM=M\text{ for every }s\in S\}.
\end{equation*}
In particular, $\mathcal{DI}:=\mathcal{D}_{R^{\times }}$ is the class of
\emph{divisible} $R$-modules. For any unexplained definitions and
terminology on domains and their modules we refer to \cite{FS2001}.

It is well known that every module over any ring has an \emph{injective
envelope} as shown by B. Eckmann and A. Schopf \cite{ES1953} (see \cite[17.9]
{Wis1991}). The dual result does not hold for the categorical dual notion of
\emph{projective covers}. Rings over which every (finitely generated)$\;$module has
a projective cover were considered first by H. Bass \cite{Bas1960} and
called (\emph{semi}-)\emph{perfect rings}. At the beginning of the current
century, L. Bican, R. El Bashir, and E. Enochs \cite{BEE2002} solved the so-called \emph{flat cover conjecture} proving that every module has a flat
cover. Recalling that the class of strongly flat modules $\mathcal{SFL}$
lies strictly between $\mathcal{FL}$ and $\mathcal{PR},$ rings over which
every (finitely generated torsion)\ module has a \emph{strongly flat cover}
were studied by S. Bazzoni and L. Salce \cite{BS2002}; such rings were
characterized as being \emph{almost }(\emph{semi-})\emph{perfect,} in the
sense that every proper homomorphic image of such rings is (semi-)perfect (see also
\cite{BS2003}). Since almost perfect rings that are not domains are perfect,
and since perfect domains are fields, the interest is restricted to almost
perfect domains (\emph{APD}'s). Although local APD's were studied earlier by
R. Smith \cite{Smi1969} under the name ``\emph{local domains with
topologically }$T$-\emph{nilpotent radical}''\emph{\ }(\emph{local
TTN-domains}), the interest in them resurfaced only recently in connection
with the revival of theory of \emph{cotorsion pairs} introduced by L. Salce
\cite{Sal1979}. Our main reference on APD's and their modules is the survey by L. Salce \cite{Sal}
(see also \cite{BS2002}, \cite{Zan2002}, \cite{BS2003}, \cite{SZ2004}, \cite {Sal2005}, \cite{Sal2006}, \cite{Zan2008}, \cite{FL2009}).

\emph{Tilting modules}\ were introduced by S. Brenner and M. Butler \cite
{BB1980} and then generalized by several authors (e.g. \cite{HR1982}, \cite
{Miy1986}, \cite{CT1995}, \cite{Wis1998}, \cite{A-HC2001}). Cotilting
modules appeared as vector space duals of tilting modules over finite
dimensional (Artin) algebras (e.g. \cite[IV.7.8.]{Hap1988}) and then
generalized in a number of papers (e.g. \cite{CDT1997}, \cite{A-HC2001},
\cite{Wis2002}, \cite{Baz2004}). A classification of (co)tilting modules
over special classes of commutative rings and domains was initiated by R.
G\"{o}bel and Trlifaj \cite{GT2000}, who classified (co)tilting Abelian
groups (assuming G\"{o}del's axiom of constructibility; a condition removed
later in \cite{BET2005}). (Co)tilting modules were classified also over
Dedekind domains by S. Bazzoni et al. \cite{BET2005} (removing set
theoretical assumptions in \cite{TW2002}), over valuation and Pr\"{u}fer
domains by L. Salce in \cite{Sal2004} and \cite{Sal2005B}, and recently over
arbitrary $1$-Gorenstein rings by J. Trlifaj and D. Posp\'{i}\v{s}il \cite
{TP2009}.

An open problem in \cite[Page 254]{GT2006} is ``\emph{Characterize all
tilting modules and classes over Matlis domains}'' ($R$ is Matlis, iff $\mathrm{%
p.d.}_{R}(Q)=1$). Recalling that APD's are Matlis domains by
\cite[Proposition 2.5]{Sal}, a natural question in this connection was
raised to the first author by L. Salce: ``\emph{Characterize all tilting
modules and classes over APD's}''. Our main result (Theorem \ref{MAIN})
provides a complete answer:

\vspace*{2mm} \textbf{MAIN THEOREM.} Let $R$ be an APD that is not a field.

\begin{enumerate}
\item  All tilting $R$-modules are $1$-tilting and represented (up to
equivalence) by
\begin{equation*}
\{T(X):=\bigcap\limits_{\frak{m}\in X}R_{\frak{m}}\bigoplus \frac{%
\bigcap\limits_{\frak{m}\in X}R_{\frak{m}}}{R}\mid X\subseteq \mathrm{Max}%
(R)\}.
\end{equation*}

\item  $\{X$\textrm{-}$\mathrm{Div}\mid X\subseteq \mathrm{Max}(R)\}$ is the
class of all tilting classes, where
\begin{equation*}
X\text{\textrm{-}}\mathrm{Div}:=\{_{R}M\mid \frak{m}M=M\text{ for every }%
\frak{m}\in X\}.
\end{equation*}

\item  If $R$ is \emph{coprimely packed}, then the set of \emph{Fuchs-Salce
tilting modules }
\begin{equation*}
\{\delta _{S}\mid S\subseteq R^{\times }\text{ is a multiplicative subset}\}
\end{equation*}
classifies all tilting $R$-modules (up to equivalence).
\end{enumerate}

This provides a partial solution to the above mentioned open problem on
Matlis domains and generalizes the classification of tilting modules over $1$%
-Gorenstein domains (which are properly contained in the class of APD's) and
Dedekind domains.

The paper is organized as follows. After this introductory section, we
collect in Section 2 some preliminaries on (semi-)perfect rings and almost
(semi-)perfect domains. In Section 3, we characterize some classes of
modules over APD's:
\begin{equation*}
\mathcal{I}=\mathcal{I}_{1},\mathcal{F}=\mathcal{F}_{1}=\mathcal{P}_{1}=%
\mathcal{P},\mathcal{IN}=\mathcal{DI}\cap \mathcal{I}_{1},\mathcal{FL}=%
\mathcal{TF}\cap \mathcal{P}_{1},\mathcal{DI}=\{M\mid \mathrm{rad}(_R M)=M\}.
\end{equation*}
Although these results are meant to serve in proving the main result (Theorem \ref{MAIN}),
we include them in a separate section since we believe they are interesting for their own.
In Section 4, we present our main results. Since $\mathcal{I}=\mathcal{I}_{1}$ and $\mathcal{P}=\mathcal{P}_{1},$ we notice
first that all (co)tilting modules over APD's are $1$-(co)tilting. Moreover,
we conclude (analogous to the case of Pr\"{u}fer domains) that all
torsion-free tilting modules over APD's are projective. In the local case,
we prove that every tilting module over a local APD is either divisible or
projective (see Theorem \ref{T=D-P}). Finally, we present in Theorem \ref
{MAIN} a complete classification of all tilting modules over APD's that are
not fields. Assuming moreover that the APD $R$ is \emph{coprimely packed}
(e.g. $R$ is a semilocal), we show that any tilting module is equivalent to
a \emph{Fuchs-Salce tilting }$R$\emph{-module} $\delta _{S}$ for some
suitable multiplicative subset $S\subseteq R^{\times }.$ If $R$ is a
coherent (whence Noetherian) APD, then the cotilting $R$-modules are
precisely the (dual) character modules of the tilting ones (see Corollary
\ref{Cot-APD}).

\section{Preliminaries}

In this section, we collect some preliminaries on (semi-)perfect rings and
almost (semi-)perfect domains.

\begin{definition}
(\cite{Bas1960}) The ring $R$ is said to be (\textbf{semi-})\textbf{perfect}%
, iff every (finitely generated) $R$-module has a projective cover.
\end{definition}

For the convention of the reader, we collect in the following lemma some of the characterizations of perfect commutative rings (e.g. \cite[Section 28]
{AF1993}, \cite[Section 43]{Wis1991}, \cite[Chapter 8]{Lam2001}, \cite[Theorem 1.1]{BS2003}):

\begin{lemma}
\emph{\label{perfect}}The following are equivalent:

\begin{enumerate}
\item  $R$ is perfect;

\item  every semisimple $R$-module has a projective cover;

\item  every flat $R$-module is \emph{(}self-\emph{)}projective;

\item  direct limits of projective $R$-modules are \emph{(}self-\emph{)}%
projective;

\item  $R$ is semilocal and every non-zero $R$-module has a maximal
submodule;

\item  $R$ is semilocal and every non-zero $R$-module contains a simple
submodule;

\item  $R$ contains no infinite set of orthogonal idempotents and every
non-zero $R$-module contains a simple submodule;

\item  $R/J(R)$ is semisimple and $J(R)$ is $T$-nilpotent;

\item  $R/J(R)$ is semisimple and $R$ is semiartinian;

\item  $R$ satisfies the DCC for principal \emph{(}finitely generated\emph{)}
ideals;

\item  Any $R$-module satisfies the DCC on its cyclic \emph{(}finitely
generated\emph{)} $R$-submodules;

\item  Any $R$-module satisfies the ACC on its cyclic $R$-submodules;

\item  $R$ is a finite direct product of local rings with $T$-nilpotent
maximal ideals;

\item  $R$ is semilocal and $R_{\frak{m}}$ is a perfect ring for every $%
\frak{m}\in \mathrm{Max}(R);$

\item  $R$ is semilocal and semiartinian.
\end{enumerate}
\end{lemma}

\begin{definition}
(\cite{BS2002}, \cite{BS2003}) $R$ is an \textbf{almost (semi-)perfect ring}%
, iff $R/I$ is (semi-)perfect for every non-zero ideal $0\neq
I\trianglelefteq R.$
\end{definition}

\begin{remark}
An almost perfect ring that is not a domain is necessarily perfect by
\cite[Proposition 1.3]{BS2003}. On the other hand, any perfect domain is a field (e.g. \cite[Corollary 1.3]{Sal}). This restricts the
interest to \emph{almost perfect domains} (\emph{APD}'s).
\end{remark}

\begin{lemma}
\label{ASPD}\emph{(\cite[Theorem 4.9]{BS2002}, \cite[Theorem IV.3.7]{FS2001})%
} The following are equivalent for an integral domain $R:$

\begin{enumerate}
\item  $R$ is almost semi-perfect;

\item  every finitely generated torsion $R$-module has a strongly flat cover;

\item  $Q/R\simeq \bigoplus\limits_{\frak{m}\in \mathrm{Max}(R)}(Q/R)_{\frak{%
m}}$ canonically;

\item  $R$ is $h$-local \emph{(}i.e. $R/I$ is semilocal for every non-zero
ideal $0\neq I\trianglelefteq R$ and $R/P$ is local for every non-zero prime
ideal $0\neq P\in \mathrm{Spec}(R)$\emph{)}.
\end{enumerate}
\end{lemma}

In the following lemma we collect several characterizations of APD's (see
\cite[Main Theorem]{Sal}, \cite{BS2002}, and \cite{BS2003}):

\begin{lemma}
\label{APD-Main}For an integral domain $R$ with $Q\neq R$ the following are
equivalent:

\begin{enumerate}
\item  $R$ is an APD;

\item  $R$ is almost semi-perfect and $R_{\frak{m}}$ is an APD for every $%
\frak{m}\in \mathrm{Max}(R);$

\item  $R$ is $h$-local and $R_{\frak{m}}$ is an APD for every $\frak{m}\in
\mathrm{Max}(R);$

\item  $R$ is $h$-local and $Q/R$ is semiartinian;

\item  $R$ is $h$-local and for every proper non-zero ideal $I\neq 0,R,$ the
$R$-module $R/I$ contains a simple $R$-submodule.

\item  every flat $R$-module is strongly flat;

\item  every $R$-module has a strongly flat cover;

\item  every weakly cotorsion $R$-module is cotorsion;

\item  every $R$-module with weak dimension at most $1$ has projective
dimension at most $1$ \emph{(}i.e. $\mathcal{F}_{1}=\mathcal{P}_{1}$\emph{)};

\item  every divisible $R$-module is weak-injective.
\end{enumerate}
\end{lemma}

\begin{remarks}
\label{rem-APD}Let $R$ be an integral domain.

\begin{enumerate}
\item  $R$ is a coherent APD if and only if $R$ is Noetherian and $1$%
-dimensional (see \cite[Propositions 4.5, 4.6]{BS2002}). Whence, Dedekind
domains are precisely the Pr\"{u}fer APD's.

\item  A valuation domain $R$ is an APD if and only if $R$ is a DVR (e.g.
\cite[Example 2.2]{Sal}).

\item  We have the following implications (e.g. \cite{FS2001}, \cite{Sal}): $%
R$ is Dedekind $\Rightarrow R$ is $1$-Gorenstein $\Rightarrow R$ is $1$%
-dimensional and Noetherian $\Rightarrow R$ is an APD $\Rightarrow R$ is a $1
$-dimensional $h$-local $\Rightarrow $ $R$ is a Matlis domain.
\end{enumerate}
\end{remarks}

The following examples illustrate that the implications above are not reversible:

\begin{examples}
\begin{enumerate}
\item  \label{ex-APD}Let $d$ be a square-free integer such that $d\equiv 1$
(mod $4$) and consider the commutative Noetherian subring
\begin{equation*}
R:=\{\frac{m}{2n+1}+\frac{m^{\prime }}{2n^{\prime }+1}\sqrt{d}\mid
m,m^{\prime },n,n^{\prime }\in \mathbb{Z}\}\subseteq \mathbb{Q}[\sqrt{d}].
\end{equation*}
By \cite[Corollary 4.5]{Smi2000}, $R$ is a $1$-Gorenstein domain that is not
Dedekind.

\item  Let $K$ be a field. Then $R=K[|t^{3},t^{5},t^{7}|]$ is a Noetherian $1
$-dimensional domain which is not $1$-Gorenstein (e.g. \cite[Ex. 18.8]{Mat}).

\item  Let $K$ be a field and $V=(K[[x]],M)$ the local domain of power
series in the indeterminate $x$ with coefficients in $K$ and with maximal
ideal $M:=xK[[x].$ Let $(D,\frak{m})$ be a local subring of $K$ and consider
the local integral domain $R:=(D+M,\frak{m}+M$). By \cite[Lemma 3.1]{BS2003}%
, $R$ is an APD if and only if $D$ is a field. Moreover, by \cite[Example
3.3]{BS2003}, if $D=F$ is a field and $K=F(X),$ then $R$ is Noetherian if
and only if $[K:F]<\infty .$ So, if $[K:F]=\infty $ then $R$ is a \emph{%
non-Noetherian} APD whence not $1$-Gorenstein.

\item  Any rank-one non-discrete valuation domain is a $1$-dimensional local
Matlis domain that is not an APD (a concrete example is \cite[Example 1.3]{Zan2008}).

\item  Any almost Dedekind domain which is not Dedekind is a $1$-dimensional Matlis domain that is not of finite character, whence not $h$-local (for a concrete example see \cite[Example III.5.5]{FS2001}).
\end{enumerate}
\end{examples}

Generalizing the so-called \emph{Prime Avoidance Theorem} (e.g. \cite[3.61]
{Sha2000}) by allowing \emph{infinite} unions of prime ideals led to the
following notions.

\begin{punto}
(\cite{RV1970}, \cite{Erd1988}) An ideal $I$ of a commutative ring $R$ is
said to be \emph{coprimely packed} (resp. \emph{compactly packed}), iff for
any set of maximal (resp. prime) $R$-ideals $\{P_{\lambda }\}_{\Lambda }$ we
have
\begin{equation}
I\subseteq \bigcup_{\lambda \in \Lambda }P_{\lambda }\Rightarrow I\subseteq
P_{\lambda _{0}}\text{ for some }\lambda _{0}\in \Lambda .  \label{packed}
\end{equation}
A class of $R$-ideals $\mathcal{E}$ said to be \emph{coprimely packed}
(resp. \emph{compactly packed}), iff every ideal in $\mathcal{E}$ is so. The
ring $R$ is said to be \emph{coprimely packed} (resp. \emph{compactly packed}%
), iff every ideal of $R$ is coprimely packed (resp. compactly packed).
\end{punto}

\begin{remark}
\label{rem-cp}By \cite[Lemma 2]{EM1994} (resp. \cite[Theorem 2.3]{BO1992}),
a ring $R$ is coprimely packed (resp. compactly packed) if and only if $%
\mathrm{Spec}(R)$ is coprimely packed (resp. compactly packed). Indeed, $1$%
-dimensional rings (e.g. APD's) are coprimely packed if and only if they are
compactly packed. By \cite{RV1970} a Dedekind domain is compactly packed
(equivalently coprimely packed) if and only if its ideal class group is
torsion (see also \cite[Theorem 1.4]{Erd1988}). Semilocal rings are
obviously coprimely packed (by the Prime Avoidance Theorem). A coprimely
packed domain $R$ is $h$-local if, for example, $R$ is $1$-dimensional by
\cite[Proposition 1.3]{Erd1988} and \cite[Theorem 3.22]{Mat1972} (see also
\cite[Theorem 3.7, EX. IV.3.3]{FS2001}) or if $Q/R$ is injective by
\cite[Theorem 9]{Bra1975}. While clearly all compactly packed rings are
coprimely packed, it had been shown in \cite{RV1970} that a Noetherian
compactly packed ring has Krull dimension at most one; thus any semilocal
Noetherian ring with Krull dimension at least $2$ is coprimely packed but
not compactly packed.
\end{remark}

\begin{example}
\label{ex-cp}Let $K$ be an algebraically closed field and $F$ a proper subfield such
that $[K:F]=\infty $ and $X$ an indeterminate. By \cite[Example 5.5]{Sal}, $%
R:=F+XK[X]$ is a non-coherent APD with $\mathrm{Max}(R)=\{XK[X]\}\cup
\{(1-aX)R\mid a\in K^{\times }\}.$ Clearly, $R$ is a coprimely packed
(compactly packed) APD that is not semilocal.
\end{example}

\section{Modules over APD's}

In this section, we characterize the injective modules, the torsion-free
modules, and the divisible modules over almost perfect domains. Moreover, we
show that over such integral domains $\mathcal{I}=\mathcal{I}_{1},$ $%
\mathcal{F}=\mathcal{F}_{1}=\mathcal{P}_{1}=\mathcal{P}.$ Throughout in this
section, $R$ is an almost perfect domain with $Q\neq R.$

Dedekind domains are characterized by the fact that every divisible module
is injective (e.g. \cite[Theorem 4.24]{Ro}, \cite[40.5]{Wis1991}). This
inspires:

\begin{proposition}
\label{inj-APD}An $R$-module $M$ is injective if and only if $M$ is
divisible and $\mathrm{i.d.}_{R}(M)\leq 1,$ i.e.
\begin{equation}
\mathcal{IN}=\mathcal{DI}\cap \mathcal{I}_{1}.  \label{IN=}
\end{equation}
\end{proposition}

\begin{Beweis}
$(\Rightarrow )$ Injective modules over any ring are divisible (e.g.
\cite[16.6]{Wis1991}).

$(\Leftarrow )$ Assume that $_{R}M$ is divisible and $\mathrm{i.d.}%
_{R}(M)\leq 1.$

\textbf{Case 1. }$(R,\frak{m})$ is \emph{local}. Let $0\neq r\in R$ be
arbitrary. By Lemma \ref{APD-Main} (5), the $R$-module $R/Rr$ contains a
simple $R$-submodule $J/Rr$ ($\simeq R/\frak{m},$ since $\mathrm{Max}(R)=\{%
\frak{m}\}$). So, we have a short exact sequence of $R$-modules
\begin{equation*}
0\rightarrow J/Rr\rightarrow R/Rr\rightarrow R/J\rightarrow 0.
\end{equation*}
Applying the contravariant functor $\mathrm{Hom}_{R}(-,M),$ we get a long
exact sequence
\begin{equation*}
\cdots \rightarrow \mathrm{Ext}_{R}^{1}(R/Rr,M)\rightarrow \mathrm{Ext}%
_{R}^{1}(J/Rr,M)\rightarrow \mathrm{Ext}_{R}^{2}(R/J,M)\rightarrow \cdots
\end{equation*}
Since $_{R}M$ is divisible, we have $\mathrm{Ext}_{R}^{1}(R/Rr,M)=0$ by
\cite[Lemma I.7.2]{FS2001}; and since $\mathrm{i.d.}_{R}(M)\leq 1,$ we have $%
\mathrm{Ext}_{R}^{2}(R/J,M)=0.$ It follows that $\mathrm{Ext}_{R}^{1}(R/%
\frak{m},M)\simeq \mathrm{Ext}_{R}^{1}(J/Rr,M)=0,$ whence $_{R}M$ is
injective by \cite[Proposition 8.1. (1)]{Sal}.

\textbf{Case 2.} $R$ is arbitrary. Let $\frak{m}\in \mathrm{Max}(R)$ be
arbitrary. Since $R$ is $h$-local, it follows by \cite[Theorem IX.7.6]
{FS2001} that localizing any injective coresolution of $R$-modules at $\frak{%
m}$ yields an injective coresolution of $R_{\frak{m}}$-modules, hence $%
\mathrm{i.d.}_{R_{\frak{m}}}(M_{\frak{m}})\leq 1.$ Since $_{R_{\frak{m}}}M_{%
\frak{m}}$ is also divisible, we conclude that $_{R_{m}}M_{\frak{m}}$ is
injective by the proof of Case 1. Since $R$ is $h$-local, we have (e.g. \cite{Mat1972}, \cite[Theorem IX.7.6]{FS2001})
\begin{equation*}
\mathrm{i.d.}_{R}(M)=\sup \{\mathrm{i.d.}_{R_{\frak{m}}}(M_{\frak{m}})\mid
\frak{m}\in \mathrm{Max}(R)\} = 0. \blacksquare
\end{equation*}
\end{Beweis}

It is well-known that for $1$-Gorenstein domains (and general $1$-Gorenstein rings), we have $\mathcal{I}= \mathcal{I}_{1} = \mathcal{F}=\mathcal{F}_{1} = \mathcal{P}=\mathcal{P}_{1}$ (e.g. \cite[9.1.10]{EJ2000}, \cite[7.1.12]{GT2000}). For the strictly larger class of APD's (see Example \ref{ex-APD} (3)), these hold partially.

\begin{proposition}
\label{I-1-P}We have
\begin{equation}
\mathcal{I}=\mathcal{I}_{1},\text{ }\mathcal{F}=\mathcal{F}_{1}=\mathcal{P}%
_{1}=\mathcal{P}.  \label{I=, P=}
\end{equation}
\end{proposition}

\begin{Beweis}
Let $R$ be an APD.

\begin{itemize}
\item  We prove, by induction, that any $R$-module $M$ with finite injective
dimension at most $n$ has injective dimension at most $1.$ If $n=0,$ we are
done. Let $n\geq 1$ and assume the statement is true for $n-1.$ Let
\begin{equation*}
0\rightarrow M\overset{f_{0}}{\longrightarrow }E_{0}\overset{f_{1}}{%
\longrightarrow }E_{1}\rightarrow \cdots \longrightarrow E_{n-2}\overset{%
f_{n-1}}{\longrightarrow }E_{n-1}\overset{f_{n}}{\longrightarrow }%
E_{n}\longrightarrow 0
\end{equation*}
be an injective coresolution of $_{R}M$ and $L:=\mathrm{Im}(f_{n-1})=\mathrm{%
Ker}(f_{n}).$ Being a homomorphic image of a divisible $R$-module, $L$ is
divisible and obviously $\mathrm{i.d.}_{R}(L)\leq 1$ whence $_{R}L$ is
injective by Proposition \ref{inj-APD}. It follows that $\mathrm{i.d.}%
_{R}(M)\leq n-1,$ whence $\mathrm{i.d.}_{R}(M)\leq 1$ by the induction
hypothesis.

\item  Let $M$ be with finite weak (flat)\ dimension at most $n.$ By
\cite[Proposition IX. 7.7]{FS2001} we have for any injective cogenerator $%
_{R}\mathbf{E}:$%
\begin{equation}
\mathrm{i.d.}_{R}(\mathrm{Hom}_{R}(M,\mathbf{E}))=\mathrm{w.d.}_{R}(M)
\label{w.d.}
\end{equation}
and we conclude that $\mathrm{w.d.}_{R}(M)\leq 1$ by the first part of the
proof.

\item  Let $_{R}M$ be with finite projective dimension at most $n.$ Since $%
\mathrm{w.d.}_{R}(M)\leq \mathrm{p.d.}_{R}(M)\leq n,$ we have $M\in \mathcal{%
F}_{1}=\mathcal{P}_{1}$ by Lemma \ref{APD-Main} (9).$\blacksquare $
\end{itemize}
\end{Beweis}

Using Proposition \ref{I-1-P} we conclude that an APD is either Dedekind or
has (weak) global dimension $\infty $. This provides new characterizations
of Dedekind domains and recovers the fact that Dedekind domains are
precisely the Pr\"{u}fer APD's.

\begin{corollary}
\label{Dedekind}An arbitrary integral domain $R$ is Dedekind if and only if $%
R$ is an APD with finite \emph{(}weak\emph{)} global dimension if and only
if $R$ is an APD with \emph{(}weak\emph{)} global dimension at most one if
and only if $R$ is a Pr\"{u}fer APD.
\end{corollary}

\begin{proposition}
\label{flat-APD}An $R$-module $M$ is flat if and only if $M$ is torsion-free
and $\mathrm{p.d.}_{R}(M)\leq 1,$ i.e.
\begin{equation}
\mathcal{FL}=\mathcal{TF}\cap \mathcal{P}_{1} = \mathcal{TF}\cap \mathcal{F}_{1}.  \label{FL=}
\end{equation}
\end{proposition}

\begin{Beweis}
$(\Rightarrow )$ Follows by the well-known fact that flat modules over
domains are torsion-free (e.g. \cite[36.7]{Wis1991}). So, we are done by $%
\mathcal{F}_{1}=\mathcal{P}_{1}$ (Lemma \ref{APD-Main} (9)).

$(\Leftarrow )$ Since $_{R}M$ is torsion-free, it embeds in a vector space
over $Q$ (e.g. \cite[Lemma 4.33]{Ro}). So, we have a short exact sequence of
$R$-modules
\begin{equation*}
0\rightarrow M\rightarrow Q^{(\Lambda )}\rightarrow Q^{(\Lambda
)}/M\rightarrow 0.
\end{equation*}
Since $_{R}Q^{(\Lambda )}$ is flat, $\mathrm{p.d.}_{R}(Q^{(\Lambda )})\leq 1$
by Lemma \ref{APD-Main} (9). It follows by \cite[Lemma VI.2.4]{FS2001} that $%
\mathrm{p.d.}_{R}(Q^{(\Lambda )}/M)<\infty ,$ whence $Q^{(\Lambda )}/M\in
\mathcal{P}_{1}=\mathcal{F}_{1}$ by Proposition \ref{I-1-P}. Consequently, $%
_{R}M$ is flat.$\blacksquare $
\end{Beweis}

\begin{punto}
(\cite{GT2006}) An $R$-module over an (arbitrary ring) $R$ is said to be%
\emph{\ }\textbf{strongly finitely presented}, iff it possesses a projective
resolution consisting of finitely generated $R$-modules. With $R$-$\mathrm{%
mod}$ we denote the class of such modules. In case $R$ is coherent, $R$-$%
\mathrm{mod}$ coincides with the class of finitely presented $R$-modules.
\end{punto}

\begin{proposition}
\label{APD-div}The following are equivalent for an $R$-module $M:$

\begin{enumerate}
\item  $_{R}M$ is divisible;

\item  $\mathrm{rad}(_{R}M)=M$ \emph{(}i.e. $M$ has no maximal $R$-submodules%
\emph{)};

\item  $\frak{m}M=M$ for every $\frak{m}\in \mathrm{Max}(R).$
\end{enumerate}
\end{proposition}

\begin{Beweis}
The result is obvious for $M=0.$ So, assume $M\neq 0.$ The equivalence $%
(1)\Leftrightarrow (3)$ is already known for APD's (e.g. L. Salce
\cite[Proposition 8.1]{Sal}).

$(1)\Rightarrow (2)$ Suppose that $M$ contains a maximal $R$-submodule $L.$ Then $M/L\simeq R/\frak{m}$
for some maximal ideal $\frak{m}\trianglelefteq R.$ Since $_{R}M$ is
divisible by assumption, it follows that $R/\frak{m}$ is also a divisible $R$-module (a contradiction).

$(2)\Rightarrow (1)$ Suppose $_{R}M$ is not divisible. Then there exists $%
0\neq r\in R$ such that $rM\neq M.$ By Lemma \ref{perfect} (5), the non-zero
$R/rR$-module $M/rM$ contains a maximal submodule $N/rM.$ Then there exists $%
\frak{m}\in \mathrm{Max}(R),$ such that
\begin{equation*}
R/\frak{m}\simeq (R/rR)/(\frak{m}/rR)\simeq (M/rM)/(N/rM)\simeq M/N.
\end{equation*}
This implies that $N\in \mathrm{Max}(_{R}M)$ (a contradiction).$\blacksquare $
\end{Beweis}

\begin{definition}
A non-empty set $\mathcal{L}$ of $R$-ideals is said to be a \textbf{%
localizing system} (or a \textbf{Gabriel topology}), iff for any ideals $%
I,J\trianglelefteq R$ we have:

(LS1) If $I\in \mathcal{L}$ and $I\subseteq J,$ then $J\in \mathcal{L};$

(LS2) If $I\in \mathcal{L}$ and $(J:_{R}r)\in \mathcal{L}$ for every $r\in I,
$ then $J\in \mathcal{L}.$
\end{definition}

\begin{definition}
Let $R$ be an integral domain and $\mathcal{E}$ be a class of $R$-ideals. We
say an $R$-module $M$ is $\mathcal{E}$-divisible, iff $IM=M$ for every $I\in
\mathcal{E}.$
\end{definition}

For any classes $\mathcal{M}$ of $R$-modules and $\mathcal{E}$ of $R$-ideals
we set
\begin{equation*}
\begin{tabular}{lll}
$\mathcal{D}(\mathcal{M})$ & $:=$ & $\{I\trianglelefteq R\mid IM=M$ for
every $M\in \mathcal{M}\};$ \\
$\mathcal{E}$\textrm{-}$\mathrm{Div}$ & $:=$ & $\{_{R}M\mid IM=M$ for every $%
I\in \mathcal{E}\}.$%
\end{tabular}
\end{equation*}
If $R$ is a domain, then $\mathcal{D}(_{R}M)$ is a localizing system by
\cite[Lemma 1.1]{Sal2005}.

\begin{lemma}
\label{F-Div}Let $R$ be an APD and $\frak{F}$ a localizing system. An $R$%
-module $M$ is $\frak{F}$-divisible if and only if $\frak{m}M=M$ for all
maximal ideals $\frak{m}$ in $\frak{F},$ i.e.
\begin{equation}
\frak{F}\text{\textrm{-}}\mathrm{Div}=(\frak{F}\cap \mathrm{Max}(R))\text{%
\textrm{-}}\mathrm{Div}.  \label{D(F)=}
\end{equation}
\end{lemma}

\begin{Beweis}
Let $M\in (\frak{F}\cap \mathrm{Max}(R))$\textrm{-}$\mathrm{Div}.$ Let $I\in
\frak{F}$ be arbitrary and set $\mathcal{M}(I):=\{\frak{m}\in \mathrm{Max}%
(R)\mid I\subseteq \frak{m}\}\subseteq \frak{F}$ by (LS1). Let $\frak{m}\in
\mathrm{Max}(R)$ be arbitrary. If $\frak{m}\in \mathcal{M}(I),$ then $\frak{m%
}_{\frak{m}}M_{\frak{m}}=(\frak{m}M)_{\frak{m}}=M_{\frak{m}}$ whence the $R_{%
\frak{m}}$-module $M_{\frak{m}}$ is divisible by Proposition \ref{APD-div},
and it follows that $(IM)_{\frak{m}}=I_{\frak{m}}M_{\frak{m}}=M_{\frak{m}}.$
On the other hand, if $\frak{m}\notin \mathcal{M}(I),$ then $I_{\frak{m}}=R_{%
\frak{m}}$ and so $(IM)_{\frak{m}}=R_{\frak{m}}M_{\frak{m}}=M_{\frak{m}}.$
Since $(IM)_{\frak{m}}=M_{\frak{m}}$ for every $\frak{m}\in \mathrm{Max}(R),$
we conclude that $IM=M$ (i.e. $M\in \frak{F}$\textrm{-}$\mathrm{Div}$).$%
\blacksquare $
\end{Beweis}

\section{Tilting and Cotilting Modules}

This section is devoted to the classification of (co)tilting modules over
APD's. For any unexplained definitions we refer to \cite{GT2006}.

For any class of $R$-modules $\mathcal{M}$ we set
\begin{equation*}
\begin{tabular}{lll}
$\mathcal{M}^{\perp _{\infty }}$ & $:=$ & $\{_{R}N\mid \mathrm{Ext}%
_{R}^{i}(M,N)=0\text{ for all }i\geq 1$ and every $M\in \mathcal{M}\};$ \\
$^{\perp _{\infty }}\mathcal{M}$ & $:=$ & $\{_{R}N\mid \mathrm{Ext}%
_{R}^{i}(N,M)=0\text{ for all }i\geq 1$ and every $M\in \mathcal{M}\};$%
\end{tabular}
\end{equation*}
Moreover, we set
\begin{equation*}
\mathcal{M}^{\perp }:=\bigcap\limits_{M\in \mathcal{M}}\mathrm{Ker}(\mathrm{%
Ext}_{R}^{1}(M,-))\text{ and }^{\perp }\mathcal{M}:=\bigcap\limits_{M\in
\mathcal{M}}\mathrm{Ker}(\mathrm{Ext}_{1}^{R}(-,M)).
\end{equation*}

\begin{punto}
For $_{R}X,$ let $\mathrm{Gen}_{n}(_{R}X)$ be the class of $R$-modules $M$
possessing an exact sequence of $R$-modules $X^{(\Lambda _{n})}\rightarrow
\cdots \rightarrow X^{(\Lambda _{1})}\rightarrow M\rightarrow 0$ (for index
sets $\Lambda _{1},\cdots ,\Lambda _{n}$). Dually, let $\mathrm{Cogen}%
_{n}(_{R}X)$ be the class of $R$-modules $M$ possessing an exact sequence of
$R$-modules $0\rightarrow M\rightarrow X^{\Lambda _{1}}\rightarrow \cdots
\rightarrow X^{\Lambda _{n}}$ (for index sets $\Lambda _{1},\cdots ,\Lambda
_{n}$). In particular, $\mathrm{Gen}(_{R}X):=\mathrm{Gen}_{1}(_{R}X)$ is the
class of $X$\emph{-generated }$R$-modules and $\mathrm{Cogen}(_{R}X):=%
\mathrm{Cogen}_{1}(_{R}X)$ is the class of $X$\emph{-cogenerated} $R$%
-modules.
\end{punto}

\begin{punto}
Let $\mathcal{A}$ and $\mathcal{B}$ be two classes of $R$-modules. Then $(%
\mathcal{A},\mathcal{B})$ is said to be a \textbf{cotorsion pair}, iff $%
\mathcal{A}=$ $^{\perp }\mathcal{B}$ and\textbf{\ }$\mathcal{B}=\mathcal{A}%
^{\perp }.$ If, moreover, $\mathrm{Ext}_{R}^{i}(A,B)=0$ for all $i\geq 1$
and $A\in \mathcal{A},$ $B\in \mathcal{B}$ we say $(\mathcal{A},\mathcal{B})$
is \textbf{hereditary}. Each class $\mathcal{M}$ of $R$-modules \emph{%
generates} a cotorsion pair $(^{\perp }(\mathcal{M}^{\perp }),\mathcal{M}%
^{\perp })$ and \emph{cogenerates} a cotorsion pair $(^{\perp }\mathcal{M}%
,(^{\perp }\mathcal{M})^{\perp }).$ For two cotorsion pairs $(\mathcal{A},%
\mathcal{B}),$ $(\mathcal{A}^{\prime },\mathcal{B}^{\prime }),$ we have $%
\mathcal{A}=\mathcal{A}^{\prime }$ if and only if $\mathcal{B}=\mathcal{B}%
^{\prime }.$
\end{punto}

\begin{punto}
An $R$-module $T$ is said to be $n$-\textbf{tilting}, iff $\mathrm{Gen}%
_{n}(_{R}T)=T^{\perp _{\infty }};$ the \textbf{induced }$n$\textbf{-tilting
class} $T^{\perp _{\infty }}$ cogenerates a \emph{hereditary cotorsion pair }%
$(^{\perp }(T^{\perp _{\infty }}),T^{\perp _{\infty }})$ with $\mathcal{A}:=$
$^{\perp }(T^{\perp _{\infty }})\subseteq \mathcal{P}_{n}$ by \cite[Lemma
5.1.8]{GT2006} (in particular, $\mathrm{p.d.}_{R}(T)\leq n$). By \cite[Lemma
6.1.2]{GT2006} (see also \cite[Theorem 3.11]{Baz2004}), $_{R}T$ is $1$%
-tilting if and only if $\mathrm{Gen}(_{R}T)=T^{\perp }.$ An $R$-module $T$
is \textbf{tilting}, iff $T$ is $n$-tilting for some $n\geq 0.$ Two tilting $%
R$-modules $T_{1},$ $T_{2}$ are said to be \textbf{equivalent} ($T_{1}\sim T_{2}$), iff
$T_{1}^{\perp _{\infty }}=T_{2}^{\perp _{\infty }}.$
\end{punto}

\begin{punto}
An $R$-module $C$ is said to be $n$-\textbf{cotilting}, iff $\mathrm{Cogen}%
_{n}(_{R}C)=$ $^{\perp _{\infty }}C;$ the \textbf{induced} $n$\textbf{%
-cotilting class }$^{\perp _{\infty }}C$ generates a \emph{hereditary
cotorsion pair} $(^{\perp _{\infty }}C,(^{\perp _{\infty }}C)^{\perp })$
with $\mathcal{B}:=(^{\perp _{\infty }}C)^{\perp }\subseteq \mathcal{I}_{n}$
by \cite[Lemma 8.1.4]{GT2006} (in particular, $\mathrm{i.d.}_{R}(C)\leq n$).
By \cite[Lemma 8.2.2]{GT2006} (see also \cite[Theorem 3.11]{Baz2004}), $_{R}C
$ is $1$-cotilting if and only if $\mathrm{Cogen}(_{R}C)=$ $^{\perp }C.$ An $%
R$-module $C$ is said to be \textbf{cotilting}, iff $C$ is $n$-cotilting for
some $n\geq 0.$ Two cotilting $R$-modules $C_{1},$ $C_{2}$ are said to be \textbf{equivalent} ($C_{1}\sim C_{2}$), iff $^{\perp _{\infty }}C_{1}=$ $^{\perp
_{\infty }}C_{2}.$
\end{punto}

\begin{remark}
Obviously, the $0$-tilting modules are precisely the projective generators,
while the $0$-cotilting modules are precisely the injective cogenerators.
\end{remark}

\begin{example}
Let $R$ be an integral domain, $S\subseteq R^{\times }$ a multiplicative
subset, and $\omega =()$ be the empty sequence. Let $F$ be the \emph{free} $R
$-module with basis
\begin{equation*}
\beta :=\{(s_{0},\cdots ,s_{n})\mid n\geq 0\text{ and }s_{j}\in S\text{ for }%
0\leq j\leq n\}\cup \{\omega \}
\end{equation*}
and $G$ the $R$-submodule of $F$ (which is in fact \emph{free}) generated by
\begin{equation*}
\{(s_{0},\cdots ,s_{n})s_{n}-(s_{0},\cdots ,s_{n-1})\mid n>0\text{ and }%
s_{j}\in S\text{ for }0\leq j\leq n\}\cup \{(s)s-\omega \}.
\end{equation*}
The $R$-module $\delta _{S}:=F/G$ is a $1$-tilting $R$-module with $\delta
_{S}^{\perp }=\mathrm{Gen}(\delta _{S})=\mathcal{D}_{S}$ as shown in \cite
{FS1992} and we call it the \textbf{Fuchs-Salce module.} It generalizes the
\textbf{Fuchs module }$\delta :=\delta _{R^{\times }}$ (introduced in \cite
{Fuc1984}), which was studied and shown to be $1$-tilting with $\delta
^{\perp }=\mathrm{Gen}(_{R}\delta )=\mathcal{DI}$ by A. Facchini in \cite
{Fac1987} and \cite{Fac1988}.
\end{example}

\begin{definition}
(\cite{GT2006}) A \textbf{Matlis localization} of the commutative ring $R$
is $S^{-1}R,$ where $S\subseteq R^{\times }$ is a multiplicative subset and $%
\mathrm{p.d.}_{R}(S^{-1}R)\leq 1.$
\end{definition}

\begin{lemma}
\label{eqv}\emph{(\cite[Proposition 5.2.24]{GT2006}, \cite[Theorem 1.1]
{A-HHT2005}) }Let $R$ be a commutative ring and $S\subseteq R^{\times }$ a
multiplicative subset.

\begin{enumerate}
\item  Let $T$ be an $n$-tilting $R$-module, $\mathcal{T}:=T^{\perp _{\infty
}}$ the induced $n$-tilting class and
\begin{equation*}
\mathcal{T}_{S}:=\{_{S^{-1}R}N\mid N\simeq S^{-1}M\text{ for some }M\in
\mathcal{T}\}.
\end{equation*}
Then $S^{-1}T$ is an $n$-tilting $S^{-1}R$-module and its induced $n$-tilting
class is
\begin{equation*}
(S^{-1}T)^{\perp _{\infty }}:=\bigcap_{i\geq 1}\mathrm{Ker}(\mathrm{Ext}%
_{S^{-1}R}^{i}(S^{-1}T,-))=\mathcal{T}_{S}=T^{\perp _{\infty }}\cap S^{-1}R%
\text{\textrm{-}}\mathrm{Mod}\text{.}
\end{equation*}
Moreover, $_{R}M\in \mathcal{T}$ if and only if $M_{\frak{m}}\in \mathcal{T}%
_{\frak{m}}$ for every $\frak{m}\in \mathrm{Max}(R).$ If $T^{\prime }$ is
another $n$-tilting $R$-module, then
\begin{equation}
T\sim T^{\prime }\Leftrightarrow T_{\frak{m}}\sim T_{\frak{m}}^{\prime }%
\text{ for all maximal ideals }\frak{m}\in \mathrm{Max}(R).  \label{TeqvT'}
\end{equation}

\item  The following are equivalent:

\begin{enumerate}
\item  $\mathrm{p.d.}_{R}(S^{-1}R)\leq 1$ \emph{(}i.e. $S^{-1}R$ is a Matlis
localization\emph{)};

\item  $T(S):=S^{-1}R\oplus \frac{S^{-1}R}{R}$ is a $1$-tilting $R$-module;

\item  $\mathrm{Gen}(_{R}S^{-1}R)=\mathcal{D}_{S}.$
\end{enumerate}

Moreover, in this case $T(S)^{\perp _{\infty }}=\mathrm{Gen}(T(S))=\mathcal{D%
}_{S}.$
\end{enumerate}
\end{lemma}

We prove now some fundamental properties of (co)tilting modules over APD's,
some of which are analogous to the case of Pr\"{u}fer domains:

\begin{proposition}
\label{tilt-1-cotilt}Let $R$ be an APD with $R\neq Q.$

\begin{enumerate}
\item  All tilting $R$-modules are $1$-tilting.

\item  The torsion-free tilting $R$-modules are precisely the projective
generators \emph{(}i.e. the $0$-tilting $R$-modules\emph{) }and are all
equivalent to $R.$

\item  Every divisible tilting $R$-modules generates $\mathcal{DI},$ whence
is equivalent to $\delta .$

\item  All localizations of $R$ are Matlis localizations. For every
multiplicative subset $S\subseteq R^{\times }$ we have a tilting $R$-module $%
T(S):=S^{-1}R\oplus S^{-1}R/R\sim \delta _{S}$ and a cotilting $R$-module $%
T(S)^{c}\sim \delta _{S}^{c}.$

\item  All cotilting $R$-modules are $1$-cotilting.

\item  The divisible cotilting $R$-modules are precisely the injective
cogenerators \emph{(}i.e. the $0$-cotilting $R$-modules\emph{)} and are
equivalent to $R^{c}:=\mathrm{Hom}_{\mathbb{Z}}(R,\mathbb{Q}/\mathbb{Z}).$
\end{enumerate}
\end{proposition}

\begin{Beweis}
\begin{enumerate}
\item  Follows directly from $\mathcal{P}=\mathcal{P}_{1}$ (\ref{I=, P=}).

\item  If $_{R}T$ is a torsion-free tilting $R$-module, then by ``1'': $T\in
\mathcal{TF}\cap \mathcal{P}_{1}\overset{\text{(\ref{FL=})}}{=}\mathcal{FL},$
whence $_{R}T$ is projective (since flat $1$-tilting modules over arbitrary
rings are projective by \cite[Corollary 2.8]{BH}). In this case, $\mathrm{Gen%
}(_{R}T)=T^{\perp }=R$\textrm{-}$\mathrm{Mod}=R^{\perp };$ consequently, $%
_{R}T$ is a projective generator and $T\sim R.$

\item  Recall that $\mathcal{F}_{1}$ generates a cotorsion pair $(\mathcal{F}%
_{1},\mathcal{WI}),$ where (by definition) $\mathcal{WI}:=\mathcal{F}%
_{1}^{\perp }$ is the class of \emph{weak-injective} $R$-modules. Notice that conditions (8) and (9) of Lemma \ref
{APD-Main} can be expressed as $(\mathcal{F}_{1},\mathcal{WI})=(\mathcal{P}%
_{1},\mathcal{DI}).$ Let $T$ be a tilting $R$-module and consider the
induced cotorsion pair $(^{\perp }(T^{\perp }),T^{\perp }).$ If $_{R}T$ is
divisible, then $T^{\perp }=\mathrm{Gen}(_{R}T)\subseteq \mathcal{DI},$
whence $\mathcal{P}_{1}=$ $^{\perp }\mathcal{DI}\subseteq $ $^{\perp
}(T^{\perp })\subseteq \mathcal{P}_{1}.$ So, $\delta ^{\perp }=\mathcal{DI}=%
\mathcal{P}_{1}^{\perp }=T^{\perp }=\mathrm{Gen}(_{R}T),$ i.e. $T$ generates
$\mathcal{DI}$ and $T\sim \delta .$

\item  For every multiplicative subset $S\subseteq R^{\times },$ the
localization $S^{-1}R$ is a flat $R$-module whence $\mathrm{p.d.}%
_{R}(S^{-1}R)\leq 1$ by Lemma \ref{APD-Main} (9). It follows by Lemma \ref
{eqv} (2)\ that $T(S):=S^{-1}R\oplus \frac{S^{-1}R}{R}$ is a tilting $R$%
-module with $T(S)^{\perp }=\mathcal{D}_{S}=\delta _{S}^{\perp },$ whence $%
T(S)\sim \delta _{S}.$ The character module of any tilting $R$-module is
cotilting by \cite[Theorem 8.1.2]{GT2006}, whence $T(S)^{c}$ is a cotilting $%
R$-module which is equivalent to $\delta _{S}^{c}$ (e.g. \cite[Theorem 8.1.13]{GT2006}%
).

\item  Follows directly from $\mathcal{I}=\mathcal{I}_{1}$ (\ref{I=, P=}).

\item  If $_{R}C$ is a divisible cotilting $R$-module, then by ``6'': $C\in
\mathcal{DI}\cap \mathcal{I}_{1}\overset{\text{(\ref{IN=})}}{=}\mathcal{IN}.$
In this case, $\mathrm{Cogen}(_{R}C)=$ $^{\perp }C=R$\textrm{-}$\mathrm{Mod}=
$ $^{\perp }R^{c};$ consequently, $_{R}C$ is an injective cogenerator and $%
C\sim R^{c}.\blacksquare $
\end{enumerate}
\end{Beweis}

The following is a key-result that will be used frequently in the sequel.

\begin{theorem}
\label{T=D-P}Let $(R,\frak{m})$ be a local APD with $R\neq Q.$ Any tilting $R
$-module is either projective or divisible. Hence, $R$ has exactly two
tilting modules $\{R,\delta \}$ \emph{(}up to equivalence\emph{)} and
exactly two tilting classes $\{R$\textrm{-}$\mathrm{Mod},$ $\mathcal{DI}\}.$
\end{theorem}

\begin{Beweis}
Let $T$ be a tilting $R$-module and assume that $_R T$ is not divisible. Then $%
T\neq 0$ and contains by Proposition \ref{APD-div} a maximal $R$-submodule $N
$ such that $T/N\simeq R/\frak{m}.$ By \cite{BS2007} all tilting modules
(over arbitrary rings) are of finite type. So, there exists $\mathcal{S}%
\subseteq \mathcal{P}_{1}\cap R$-$\mathrm{mod}$ such that $R/\frak{m}\in
\mathrm{Gen}(_{R}T)=T^{\perp }=\mathcal{S}^{\perp }.$ Let $M\in \mathcal{S}$
be arbitrary, so that $\mathrm{Ext}_{R}^{1}(M,R/\frak{m})=0.$ Since the field $R/\frak{m}$
is indeed injective as a module over itself, it follows (e.g. \cite[Page 34 (6)]{FS2001}) that
\begin{equation*}
\begin{tabular}{lll}
$\mathrm{Tor}_{1}^{R}(R/\frak{m},M)$ & $\simeq $ & $\mathrm{Tor}_{1}^{R}(%
\mathrm{Hom}_{R/\frak{m}}(R/\frak{m},R/\frak{m}),M)$ \\
& $\simeq $ & $\mathrm{Hom}_{R/\frak{m}}(\mathrm{Ext}_{R}^{1}(M,R/\frak{m}%
),R/\frak{m})=0.$%
\end{tabular}
\end{equation*}
By \cite[II.3.2.Corollary 2]{B}, $_{R}M$ is projective (being finitely
presented and flat). So, $\mathcal{S}\subseteq \mathcal{PR},$ whence $_{R}T$
is projective.$\blacksquare $
\end{Beweis}

Recall (from \cite{Ham1971}) that an $R$-submodule $M$ of an $R$-module $N$
is said to be a \textbf{restriction submodule}, iff $M_{\frak{m}}=N_{\frak{m}%
}$ or $M_{\frak{m}}=0$ for every $\frak{m}\in \mathrm{Max}(R).$ For any
subset $X\subseteq \mathrm{Max}(R),$ we set
\begin{equation*}
R_{(X)}:=\bigcap\limits_{\frak{m}\in X}R_{\frak{m}}\text{ (}:=Q,\text{ if }%
X=\varnothing \text{) .}
\end{equation*}

\begin{lemma}
\label{loc-int}Let $R\neq Q,$ $X\subseteq \mathrm{Max}(R),$ $X^{\prime }:=%
\mathrm{Max}(R)\backslash X$ and consider
\begin{equation*}
M_{1}:=\frac{R_{(X)}}{R}\text{ and }M_{2}:=\frac{R_{(X^{\prime })}}{R}.
\end{equation*}

\begin{enumerate}
\item  If $R$ is an $h$-local domain, then $M_{1},M_{2}\subseteq \frac{Q}{R}$
are restriction $R$-submodules and
\begin{equation}
\frac{Q}{R}=M_{1}\oplus M_{2}=\frac{R_{(X)}}{R}\oplus \frac{R_{(X^{\prime })}%
}{R}.  \label{Q/R}
\end{equation}

\item  If $R$ is a $1$-dimensional $h$-local domain, then
\begin{equation*}
T(X):=R_{(X)}\bigoplus \frac{R_{(X)}}{R}\text{ \ \emph{(}}=Q\oplus \frac{Q}{R%
},\text{ if }X=\varnothing \text{\emph{)}}
\end{equation*}
is a $1$-tilting $R$-module.
\end{enumerate}
\end{lemma}

\begin{Beweis}
Recall first that if $\frak{m},\frak{m}^{\prime }\in \mathrm{Max}(R)$ are
such that $\frak{m}\neq \frak{m}^{\prime },$ then we have by \cite[Theorem
3.19]{Mat1972} (see also \cite[IV.3.2]{FS2001}):
\begin{equation}
R_{\frak{m}}\otimes _{R}R_{\frak{m}^{\prime }}\simeq (R_{\frak{m}})_{\frak{m}%
^{\prime }}=Q.  \label{R_m_m'}
\end{equation}
Moreover, if $\{R_{\lambda }\}_{\Lambda }$ is a class of $R$-submodules of $Q
$ with $\bigcap_{\lambda \in \Lambda }R_{\lambda }\neq 0,$ then it follows
from \cite[IV.3.10]{FS2001} that
\begin{equation}
(\bigcap_{\lambda \in \Lambda }R_{\lambda })_{\frak{m}}=\bigcap_{\lambda \in
\Lambda }(R_{\lambda })_{\frak{m}}\text{ for every }\frak{m}\in \mathrm{Max}%
(R).  \label{int-loc}
\end{equation}

\begin{enumerate}
\item  Clearly $M_{1}\cap M_{2}=0.$ Let $\frak{m}^{\prime }\in \mathrm{Max}%
(R)$ be arbitrary. Then
\begin{equation*}
(M_{1})_{\frak{m}^{\prime }}=\frac{(R_{(X)})_{\frak{m}^{\prime }}}{R_{\frak{m%
}^{\prime }}}\overset{\text{(\ref{int-loc})}}{=}\frac{\bigcap\limits_{\frak{m%
}\in X}(R_{\frak{m}})_{\frak{m}^{\prime }}}{R_{\frak{m}^{\prime }}}\overset{%
\text{(\ref{R_m_m'})}}{=}\left\{
\begin{tabular}{lll}
$0,$ &  & $\frak{m}^{\prime }\in X$ \\
&  &  \\
$\frac{Q}{R_{\frak{m}^{\prime }}},$ &  & $\frak{m}^{\prime }\notin X$%
\end{tabular}
\right. .
\end{equation*}
\newline
Similarly,
\begin{equation*}
(M_{2})_{\frak{m}^{\prime }}=\left\{
\begin{tabular}{lll}
$\frac{Q}{R_{\frak{m}^{\prime }}},$ &  & $\frak{m}^{\prime }\in X$ \\
&  &  \\
$0,$ &  & $\frak{m}^{\prime }\notin X$%
\end{tabular}
\right. .
\end{equation*}
So, $M_{1},M_{2}\subseteq \frac{Q}{R}$ are restriction $R$-submodules.
Moreover, we have $(M_{1}\oplus M_{2})_{\frak{m}^{\prime }}=(M_{1})_{\frak{m}%
^{\prime }}\oplus (M_{2})_{\frak{m}^{\prime }}=\frac{Q}{R_{\frak{m}^{\prime
}}}=(\frac{Q}{R})_{\frak{m}^{\prime }}$ for all $\frak{m}^{\prime }\in
\mathrm{Max}(R),$ and so $\frac{Q}{R}=M_{1}\oplus M_{2}.$

\item  Notice first that a $1$-dimensional $h$-local domain is a Matlis
domain (in fact $\mathrm{p.d.}_{R}(Q)=\mathrm{p.d.}_{R}(\frac{Q}{R})=1$ as
shown in \cite[Lemma 2.4]{Sal}). For any $X\subseteq \mathrm{Max}(R),$ we
have $\frac{Q}{R}\overset{\text{(\ref{Q/R})}}{=}\frac{R_{(X)}}{R}\oplus
\frac{R_{(X^{\prime })}}{R}$ and so $T(X)$ is a $1$-tilting $R$-module by
\cite[Theorem 8.2]{A-HHT2005}.$\blacksquare $
\end{enumerate}
\end{Beweis}

\begin{remark}
Although we proved (\ref{Q/R}) for general $h$-local domains, we point out
here that it can be obtained for an \emph{APD} $R$ by applying
\cite[Theorem 3.10]{A-HHT2005} to $M_{1}:=\frac{R_{(X)}}{R}.$ Then $X_{1}:=%
\mathrm{Supp}(M_{1})=\mathrm{Max}(R)\backslash X$ and $X_{2}:=\mathrm{Supp}%
(Q/R)\backslash X_{1}=X.$ Consider the embedding $\varphi :\frac{Q}{R}%
\rightarrow \prod\limits_{\frak{m}\in \mathrm{Max}(R)}(\frac{Q}{R})_{\frak{m}%
}.$ Since $R$ is $h$-local, it follows by \cite[Theorem IV.3.7]{FS2001} (3)
that $M_{1}\simeq \bigoplus\limits_{\frak{m}\notin \mathrm{Max}(R)}(M_{1})_{%
\frak{m}}=\bigoplus\limits_{\frak{m}\in X}\frac{Q}{R_{\frak{m}}}.$ So, $%
M_{2}:=\varphi ^{-1}(\prod\limits_{\frak{m}\in X}(\frac{Q}{R})_{\frak{m}})=%
\frac{R_{(X^{\prime })}}{R}.$ Notice that $\mathrm{w.d.}_{R}(\frac{Q}{R_{(X)}%
})\leq 1$ and so $\mathrm{p.d.}_{R}(\frac{Q}{R_{(X)}})\leq 1$ by Lemma \ref
{APD-Main} (9). The equality (\ref{Q/R}) follows now by \cite[Theorem 3.10]
{A-HHT2005}.
\end{remark}

\begin{lemma}
\label{1-dim-h-local}Let $R$ be an APD with $R\neq Q.$ If $T$ is a tilting $R
$-module, then
\begin{equation}
T^{\perp _{\infty }}=\mathrm{Gen}(_{R}T)=\mathcal{D}(_{R}T)\text{-}\mathrm{%
Div}.  \label{D(D(T))}
\end{equation}
\end{lemma}

\begin{Beweis}
Clearly $\mathrm{Gen}(_{R}T)\subseteq \mathcal{D}(T)$\textrm{-}$\mathrm{Div}.
$ Let $M\in \mathcal{D}(T)$\textrm{-}$\mathrm{Div},$ $\frak{m}\in \mathrm{Max%
}(R)$ be arbitrary and consider the tilting $R_{\frak{m}}$-module $T_{\frak{m%
}}.$ By Theorem \ref{T=D-P}, $_{R_{\frak{m}}}T_{\frak{m}}$ is either
divisible or projective. If $\frak{m}\in \mathcal{D}(T),$ then $T_{\frak{m}}$
is divisible and generates all divisible $R_{\frak{m}}$-modules by
Proposition \ref{tilt-1-cotilt} (3). Moreover, $\frak{m}_{\frak{m}}M_{\frak{m%
}}=({\frak{m}}M)_{\frak{m}}=M_{\frak{m}}$ and it follows by Proposition \ref
{APD-div} that $M_{\frak{m}}$ is a divisible $R_{\frak{m}}$-module, whence $%
M_{\frak{m}}\in \mathrm{Gen}(_{R_{\frak{m}}}T_{\frak{m}}).$ On the other
hand, if $\frak{m}\notin \mathcal{D}(T)$ then $T_{\frak{m}}$ is a projective
$R_{\frak{m}}$-module whence a generator in $R_{\frak{m}}$\textrm{-}$\mathrm{%
Mod}$ by Proposition \ref{tilt-1-cotilt} (2). In either cases $M_{\frak{m}%
}\in \mathrm{Gen}(_{R_{\frak{m}}}T_{\frak{m}})=T_{\frak{m}}^{\perp _{\infty
}}\;$for every $\frak{m}\in \mathrm{Max}(R),$ whence $M\in T^{\perp _{\infty
}}=\mathrm{Gen}(_{R}T)$ by Lemma \ref{eqv} (1).$\blacksquare $
\end{Beweis}

\begin{theorem}
\label{MAIN}Let $R$ be an APD with $R\neq Q.$

\begin{enumerate}
\item  The set
\begin{equation*}
\{T(X)\mid X\subseteq \mathrm{Max}(R)\}
\end{equation*}
is a representative set \emph{(}up to equivalence\emph{) }of all tilting $R$%
-modules.

\item  There is a bijective correspondence between the set of all tilting
torsion classes of $R$-modules and the power set of the maximal spectrum $%
\frak{B}(\mathrm{Max}(R)).$ The correspondence is given by the mutually
inverse assignments:
\begin{equation*}
\begin{tabular}{lllll}
& $\mathcal{T}\mapsto $ & $\mathcal{DM}(\mathcal{T})$ & $:=$ & $\{\frak{m}%
\in \mathrm{Max}(R)\mid \frak{m}M=M$ for every $M\in \mathcal{T}\};$ \\
and &  &  &  &  \\
& $X\mapsto $ & $X$\textrm{-}$\mathrm{Div}$ & $:=$ & $\{_{R}M\mid \frak{m}M=M
$ for every $\frak{m}\in X\}.$%
\end{tabular}
\end{equation*}

\item  If $R$ is coprimely packed, then the class of Fuchs-Salce tilting
modules
\begin{equation*}
\{\delta _{S}\mid S\subseteq R^{\times }\text{ is a multiplicative subset}\}
\end{equation*}
classifies all tilting $R$-modules \emph{(}up to equivalence\emph{)}.
\end{enumerate}
\end{theorem}

\begin{Beweis}
\begin{enumerate}
\item  Let $T$ be a tilting $R$-module and set
\begin{equation*}
\begin{tabular}{lll}
$\Omega _{1}$ & $:=$ & $\{\frak{m}\in \mathrm{M}\mathrm{ax}(R)\mid T_{\frak{m%
}}$ is a divisible $R_{\frak{m}}$-module$\};$ \\
$\Omega _{2}$ & $:=$ & $\{\frak{m}\in \mathrm{M}\mathrm{ax}(R)\mid T_{\frak{m%
}}$ is a projective $R_{\frak{m}}$-module$\}.$%
\end{tabular}
\end{equation*}
Notice first that $\mathrm{M}\mathrm{ax}(R)=\Omega _{1}\cup \Omega _{2}$ by
Theorem \ref{T=D-P} (a disjoint union by applying Proposition \ref
{tilt-1-cotilt} (2) $\&\;$(3) to the ring $R_{\frak{m}}$).

\textbf{Claim}: $T\sim T(\Omega _{2}).$ One can show (as in the proof of
Lemma \ref{loc-int}), that if $\frak{m}\in \mathrm{Max}(R)$ then
\begin{equation*}
T(\Omega _{2})_{\frak{m}}=\left\{
\begin{tabular}{lll}
$Q\oplus \frac{Q}{R_{\frak{m}}},$ &  & $\frak{m}\in \Omega _{1}$ \\
&  &  \\
$R_{\frak{m}},$ &  & $\frak{m}\in \Omega _{2}$%
\end{tabular}
\right. .
\end{equation*}
So, $T_{\frak{m}}\sim T(\Omega _{2})_{\frak{m}}$ for every $\frak{m}\in
\mathrm{Max}(R)$ whence $T\sim T(\Omega _{2})$ by (\ref{TeqvT'}).

\item  Let $\mathcal{T}=T^{\perp _{\infty }}$ be a tilting torsion class for
some tilting $R$-module $T.$ Then
\begin{equation*}
\mathcal{DM}(\mathcal{T})\text{\textrm{-}}\mathrm{Div}=\mathcal{DM}(T)\text{%
\textrm{-}}\mathrm{Div}\overset{\text{(\ref{D(F)=})}}{=}\mathcal{D}(T)\text{%
\textrm{-}}\mathrm{Div}\overset{\text{(\ref{D(D(T))})}}{=}\mathrm{Gen}%
(_{R}T)=T^{\perp _{\infty }}=\mathcal{T}.
\end{equation*}
On the other hand, let $X\subseteq \mathrm{Max}(R),$ $\overline{X}:=\mathrm{%
Max}(R)\backslash X,$ and $T^{\prime }:=T(\overline{X}).$ Then clearly $%
\mathcal{DM}(T^{\prime })=X$ and so
\begin{equation*}
\mathcal{DM}(X\text{\textrm{-}}\mathrm{Div})=\mathcal{DM}(\mathcal{DM}%
(T^{\prime })\text{\textrm{-}}\mathrm{Div})=\mathcal{DM}(T^{\prime })=X.
\end{equation*}

\item  Let $R$ be compactly packed. Let $\Omega _{1}$ and $\Omega _{2}$ be
as in ``1''.

\textbf{Case 1.} $\mathrm{Max}(R)=\Omega _{1}$ (i.e. $T_{\frak{m}}$ is a
divisible $R_{\frak{m}}$-module for all $\frak{m}\in \mathrm{Max}(R)$). In
this case, $_{R}T$ is divisible whence $T\sim Q\oplus Q/R$ and we can take $%
S=R^{\times }.$

\textbf{Case 2.} $\mathrm{Max}(R)=\Omega _{2}$ (i.e. $T_{\frak{m}}$ is a
projective $R_{\frak{m}}$-module for all $\frak{m}\in \mathrm{Max}(R)$). In
this case, $_{R}T$ is projective whence $T\sim R$ and we can take $S=\{1\}.$

\textbf{Case 3.} $\mathrm{Max}(R)\neq \Omega _{1}$ and $\mathrm{Max}(R)\neq
\Omega _{2}.$ Let
\begin{equation*}
S:=R\backslash \bigcup\limits_{\frak{m}\in \Omega _{2}}\frak{m}\text{ and }%
T(S):=S^{-1}R\oplus S^{-1}R/R.
\end{equation*}
Let $\frak{m}\in \Omega _{2},$ so that $T_{\frak{m}}$ is projective and $%
S\subseteq R\backslash \frak{m}.$ Then $(S^{-1}R)_{\frak{m}}=R_{\frak{m}}.$
Therefore $(T(S))_{\frak{m}}=(S^{-1}R)_{\frak{m}}\oplus (S^{-1}R/R)_{\frak{m}%
}=R_{\frak{m}}$ is equivalent to the projective $R_{\frak{m}}$-module $T_{%
\frak{m}}.$ On the other hand, let $\frak{m}\in \Omega _{1}$ so that $T_{%
\frak{m}}$ is a divisible $R_{\frak{m}}$-module. Then $\frak{m}\cap S\neq
\varnothing $ (otherwise $\frak{m}\subseteq \bigcup\limits_{\frak{m}\in
\Omega _{2}}\frak{m}$ and so $\frak{m}\in \Omega _{2}$ since $R$ is
coprimely packed; a contradiction since $\Omega _{1}\cap \Omega
_{2}=\varnothing $). Let $\widetilde{s}\in S\cap \frak{m}.$ Clearly $%
\widetilde{s}(S^{-1}R)_{\frak{m}}=(S^{-1}R)_{\frak{m}},$ whence $(S^{-1}R)_{%
\frak{m}}$ is a divisible $R_{\frak{m}}$-module by Proposition \ref{APD-div}%
. It follows that $(T(S))_{\frak{m}}=(S^{-1}R)_{\frak{m}}\oplus (S^{-1}R)_{%
\frak{m}}/R_{\frak{m}}$ is a divisible $R_{\frak{m}}$-module, whence $T(S)_{%
\frak{m}}\sim T_{\frak{m}}$ as $R_{\frak{m}}$-modules by Proposition \ref
{tilt-1-cotilt} (3) (applied to the ring $R_{\frak{m}}$). Since $T_{\frak{m}%
}\sim T(S)_{\frak{m}}$ for all $\frak{m}\in \mathrm{Max}(R),$ we conclude
that $T\sim T(S)$ by (\ref{TeqvT'}).$\blacksquare $
\end{enumerate}
\end{Beweis}

\begin{remark}
\label{Bass}Let $R$ be a $1$-Gorenstein ring and $_{R}T$ be a tilting $R$%
-module. By \cite{TP2009} there exists $X\subseteq \mathbf{P}_{1}$ (the set
of prime ideals of height $1$) and some (\emph{unique}) $R$-module $R_{X},$
satisfying $R\subseteq R_{X}\subseteq Q$ and fitting in an exact sequence
\begin{equation*}
0\rightarrow R\rightarrow R_{X}\rightarrow \bigoplus_{\frak{m}\in X}\mathrm{E%
}(R/\frak{m})\rightarrow 0,
\end{equation*}
such that $T$ is equivalent to the so-called \textbf{Bass tilting module }$%
B(X):=R_{X}\oplus \bigoplus_{\frak{m}\in X}\mathrm{E}(R/\frak{m}).$\textbf{\
}Let $\frak{m}\in \mathrm{Max}(R)$ be arbitrary. By the proof of
\cite[Theorem 0.1]{TP2009}, the $R_{\frak{m}}$-module $B(X)_{\frak{m}}$ is
injective, whence divisible, if $\frak{m}\in X$ and projective if $\frak{m}%
\notin X.$ If $R$ is a $1$-Gorenstein domain (whence an APD), the same holds
for the $R_{\frak{m}}$-module $T(X^{\prime })_{\frak{m}},$ where $X^{\prime
}:=\mathrm{Max}(R)\backslash X.$ It follows that, in this case, $B(X)\sim
T(X^{\prime })$ by (\ref{TeqvT'}) and so $T\sim T(X^{\prime }).\blacksquare $
\end{remark}

A direct application of Theorem \ref{MAIN},
and \cite[Theorem 8.2.8]{GT2006} yields

\begin{corollary}
\label{Cot-APD}Let $R$ be a coherent \emph{(}Noetherian\emph{)} APD.

\begin{enumerate}
\item  All cotilting $R$-modules are of cofinite type and $\{T(X)^{c}\mid
X\subseteq \mathrm{Max}(R)\}$ is a representative set \emph{(}up to
equivalence\emph{) }of all cotilting $R$-modules.

\item  If $R$ is coprimely packed, then $\{\delta _{S}^{c}\mid S\subseteq
R^{\times }$ is a multiplicative subset$\}$ classifies all cotilting $R$%
-modules \emph{(}up to equivalence\emph{)}.
\end{enumerate}
\end{corollary}

\end{document}